\count100= 1
\count101= 24
\magnification\magstep1

\nopagenumbers

\def\abs#1{\bigskip\bigskip\medskip
    {\narrower \baselineskip 9pt \fab \noindent {\bf Abstract.\ }
    \textfont0 = \fab
    \scriptfont0 = \sfab
    \scriptscriptfont0 = \fs
    \textfont1 = \tab
    \scriptfont1 = \sab
    \scriptscriptfont1 = \ssab
    \textfont2 = \abst
    \scriptfont2 = \sabs
    \scriptscriptfont2 = \ssabs
    \let \it = \itab
    \let \sl = \slab
#1 \bigskip}\medskip}

  % ACKNOWLEDGMENT

\def\Acknowledgements{\goodbreak\medskip\noindent{\bf
   Acknowledgements.\ }}  % ACKNOWLEDGEMENT

\def\Address{\goodbreak\bigskip\obeylines}
\def\ats{\kern .03em @ \kern .02em}

\def\aut{\bigskip \bigskip
        \centerline {\fb\author} }

\def\eop{{ \vrule height7pt width7pt depth0pt}\par\bigskip}

\headline{\ifodd\pageno\rightheadline \else\leftheadline\fi}

\def\leftheadline{\ifnum\pageno=\count100 \nullx%
  \else\rm\folio\hfil\it\authead\fi}%

\def\nullx{\hfill}

\def\ref{\global\advance\refnum by 1 \item{\the\refnum .}}
  \newcount\refnum \refnum = 0

\def\References{\goodbreak\bigskip\centerline{\bf References}\bigskip
   \frenchspacing}

\def\rightheadline{\ifnum\pageno=\count100 \nullx%
  \else\it\chptitle\hfil\rm\folio\fi}
\def\RR{\mathop{{\rm I}\kern-.2em{\rm R}}\nolimits}

 % DEFINE SECTION HEADING

\def\tit#1{\centerline {\fa #1}}

\def\titexp#1#2{\hbox{{\fa #1} \kern-.25em \raise .90ex \hbox{\fb #2}}\/}
\def\titsub#1#2{\hbox{{\fa #1} \kern-.25em \lower .60ex \hbox{\fb #2}}\/}

\newcount\eqnum
\newcount\figurenumb\figurenumb=0
\newcount\proclaimnumb\proclaimnumb=0
\newcount\refnumb
\newcount\writeq

\font\fa=cmr17
\font\fb=cmr12
\font\fab=cmr9
\font\sfab=cmr7
\font\fs=cmr6

\font\slab=cmsl9
\font\tab= cmmi9
\font\sab= cmmi7
\font\ssab= cmmi6
\font\abst= cmsy9
\font\sabs= cmsy7
\font\ssabs= cmsy6
\font\itab= cmti9
\chardef\sz="19

\null\vskip 18pt
\pageno=\count100
\count102=\count100
\advance\count102 by -1
\advance\count102 by \count101

\tit{Strictly Hermitian Positive Definite Functions}

\def\chptitle{Strictly Hermitian Functions}

\def\alp{\alpha}

\def\eps{\varepsilon}

\def\tet{\theta}
\def\lam{\lambda}

\def\calA{{\cal A}}
\def\calB{{\cal B}}

\def\calH{{\cal H}}
\def\calF{{\cal F}}
\def\bfa{{\bf a}}
\def\bfb{{\bf b}}
\def\bfc{{\bf c}}
\def\bfd{{\bf d}}
\def\bfx{{\bf x}}
\def\bfy{{\bf y}}

\def\bfw{{\bf w}}
\def\bfz{{\bf z}}
\def\bfv{{\bf v}}
\def\oz{{\overline z}}
\def\oc{{\overline c}}
\def\ow{{\overline w}}
\def\ob{{\overline b}}
\def\oa{{\overline a}}
\def\oC{{\overline C}}
\def\oV{{\overline V}}
\def\oU{{\overline U}}
\def\CC{{\rlap {\raise 0.4ex \hbox{$\scriptscriptstyle |$}}\hskip -0.13em C}}
\def\RR{{I\!\!R}}
\def\ZZ{{Z\!\!\! Z}}
\def\NN{{I\!\!N}}
\def\\{{\backslash}}

\def\union{\bigcup}
\def\nek{,\ldots,}
\def\span{{\rm span}}

\def\mod{{\rm mod}}
\def\rank{{\rm rank}}

\def\tilU{{\widetilde U}}
\def\sqr#1#2{{\vcenter{\hrule height.#2pt\hbox{\vrule
width.#2pt height#1pt \kern#1pt \vrule width.#2pt}\hrule
height.#2pt}}}

\def\author{Allan Pinkus}
\def\authead{Allan Pinkus}

\aut
\abs {Let $H$ be any complex inner product space with inner product
$<\cdot,\cdot>$. We say that $f:\CC\to\CC$ is Hermitian positive
definite on $H$ if the matrix
$$\left(f(<\bfz^r,\bfz^s>)\right)_{r,s=1}^n\eqno(*)\phantom{12341}$$
is Hermitian positive definite for all choice of $\bfz^1\nek \bfz^n$
in $H$, all $n$. It is strictly Hermitian positive definite if the
matrix ({\lower 0.2pt\hbox{$*$}})
is also non-singular for any choice of distinct
$\bfz^1\nek \bfz^n$ in $H$. In this article we prove that if $\dim
H\ge 3$, then $f$ is Hermitian positive definite on $H$ if and only if
$$f(z) = \sum_{k,\ell =0}^\infty b_{k,\ell} z^k \oz^\ell
\eqno(**)\phantom{12341}$$
where $b_{k,\ell}\ge 0$, all $k,\ell$ in $\ZZ_+$, and the series
converges for all $z$ in $\CC$. We also prove that $f$ of the form
({\lower 0.2pt\hbox{$**$}}) is  strictly Hermitian
positive definite on any $H$ if and only if the set
$$J=\{ (k,\ell):\, b_{k,\ell}> 0\}$$
is such that $(0,0)\in J$, and
every arithmetic sequence in $\ZZ$ intersects the values
$\{k-\ell\,:\, (k,\ell)\in J\}$ an infinite number of times.}

\null

\bigskip\noindent
\centerline{\bf \S 1. Introduction}

\medskip\noindent
A function (or kernel) $K$ mapping the product space $Z\times Z$
into $\CC$ is termed {\sl Hermitian positive definite} if
$$\sum_{r,s=1}^n c_r K(\bfz^r,\bfz^s) \oc_s\ge 0\,,\eqno(1.1)$$
for every choice of $\bfz^1\nek\bfz^n \in Z
$, $c_1\nek c_n\in
\CC$, and all $n\in \NN$. We say that the function $K$ is {\sl
positive definite} if $K:X\times X \to \RR$ and $$\sum_{r,s=1}^n
c_r K(\bfx^r,\bfx^s) c_s\ge 0\,,\eqno(1.2)$$ for every choice of
$\bfx^1\nek\bfx^n \in X$, $c_1\nek c_n\in \RR$, and all $n\in
\NN$. These functions seem to have been first considered by Mercer
[20] in connection with integral equations. We use the term {\sl
strict} if strict inequalities occur in (1.1) and (1.2) for every
choice of distinct $\bfz^1\nek \bfz^n\in Z$ or $\bfx^1\nek
\bfx^n\in X$, as appropriate, and nonzero $c_1\nek c_n$.

Hermitian positive  definite and positive definite functions have
been much studied in various contexts and guises. One of the first
characterizations of sets of Hermitian positive  definite
functions is Bochner's Theorem, see e.g. Bochner [4], Chung [9,
Sect.~6.5], and Kawata [16, Chap.~10]. The function $$K(x,y) =
f(x-y)\,,$$ where  $X=\RR$ and $f$ is continuous at $0$, is
Hermitian positive definite if and only if $f$ is of the form
$$f(t) =\int_{-\infty}^\infty e^{itx} d\mu(x)$$ where $d\mu$ is a
finite, nonnegative measure on $\RR$. (The essentially same
result holds for $X=\RR^n$, see Bochner [4].)

In the theory of radial basis functions there has been much
interest in characterizing positive definite functions of the form
$$K(\bfx,\bfy) = f(\|\bfx-\bfy\|)$$ where $\|\cdot\|$ is some
norm. There is a large literature connected with such problems.
Schoenberg [26] characterized such functions where $\|\cdot\|$ is
the Euclidean norm on $X=\RR^n$. For the analogous problem with
the $\ell_p$ norm $\|\cdot\|_p$ on $\RR^n$, see e.g. Gneiting [12]
and references therein. Schoenberg [27] also characterized
positive definite functions of the form $$K(\bfx,\bfy) =
f(\|\bfx-\bfy\|)$$ where $\|\cdot\|$ is again the Euclidean norm
on $\RR^n$, but $X=S^{n-1}$ is the unit sphere in $\RR^n$, $n\ge
2$. As there is a simple 1-1 correspondence between
$\|\bfx-\bfy\|$ and the standard inner product $<\bfx, \bfy>$ for
all $\bfx,\bfy\in S^{n-1}$, it is more convenient to consider
$$K(\bfx, \bfy) = g(<\bfx, \bfy>)\,.$$ Schoenberg proved that for
$g\in C[-1,1]$ the kernel $K$ is positive definite if and only if
$g$ is of the form $$g(t) = \sum_{r=0}^\infty a_r P_r^\lam
(t)\eqno(1.3)$$ where $a_r\ge 0$ for all $r$, $\sum_{r=0}^\infty
a_r P_r^\lam(1)
<\infty$, and the $P_r^\lam$ are the Gegenbauer (ultraspherical)
polynomials with $\lam = (n-2)/2$. As these kernels are often used
for interpolation, there has been much effort
put into determining exact conditions for when $g(<\cdot, \cdot>)$
is strictly positive definite. It was recently proven by Chen,
Menegatto, Sun [7] that for $n\ge 3$ a function of the form (1.3)
is strictly positive definite if and only if the set $$\{ r\,:\,
a_r>0\}$$ contains an infinite number of even and  an infinite
number of odd integers.

Hermitian positive definite and positive definite functions also
arise in the study of reproducing kernels, see e.g. Aronszajn [1],
Meschkowski [21], and Donoghue [10, Chap.~X]. Each reproducing
kernel is a (Hermitian) positive definite function and vice versa.
Strict (Hermitian) positive definiteness is also desired when
considering the reproducing kernel since it is equivalent to the
property of the linear independence of point functionals in
the associated reproducing kernel space.

As a corollary (quite literally a footnote) in Schoenberg [27] is
the result that if for every positive definite matrix
$(a_{rs})_{r,s=1}^n$, all $n\in \NN$, the function $f:\RR\to\RR$
is such that $$\left(f(a_{rs})\right)_{r,s=1}^n$$ is also positive
definite, then $f$ is necessarily of the form $$f(t)
=\sum_{k=0}^\infty b_k t^k\eqno(1.4)$$ where the $b_k\ge 0$ for all
$k$, and the series converges for all $t\in \RR$. That is, $f$ is
real entire and absolutely monotone on $\RR_+$. The converse
direction is a simple consequence of the fact that the class of
positive definite functions is a positive cone and the Schur Product Theorem, and is
an exercise in P\'olya, Szeg\"o [24, p.~101]. For other approaches
and generalizations see Rudin [25], Christensen, Ressel [8], Berg,
Christensen, Ressel [2, p.~59] and FitzGerald, Micchelli, Pinkus
[11]. The same problem for matrices of a fixed size, i.e., all
$n\times n$ matrices for a fixed $n$, seems much more difficult,
see e.~g. Horn [14].

If $H$ is a real inner product space of dimension $m$, then the
set of matrices $(<\bfx^r, \bfx^s>)_{r,s=1}^n$, all $n$, obtained
by choosing arbitrary $\bfx^r$ in $H$ is exactly the set of all
positive definite matrices of rank at most $m$. Thus if $\dim
H=\infty$, then (1.4) provides the exact characterization of all functions $f$
such that $$f(<\bfx,\bfy>)\eqno(1.5)$$ is positive definite. In
Lu, Sun [19] this result was significantly improved. They proved
that if $\dim H \ge 2$, then $f:\RR\to\RR$ is positive definite on
$H$ (in the sense of (1.5)) if and only if $f$ is of the form
(1.4). If $\dim H=1$ then this result is not valid. In Pinkus [23]
it was proven that assuming that $f$ has the form (1.4), then $f$
(as given in (1.5)) is strictly positive definite on $H$ if
and only if the set $$\{k\,:\, b_k>0\}$$ contains the index $0$
plus an infinite number of even integers and an infinite number of
odd integers.

In this paper we will generalize both these results to $\CC$. In
Herz [13], see also  Berg, Christensen, Ressel [2, p.~171], and
FitzGerald, Micchelli, Pinkus [11], it was proven that if
$f:\CC\to\CC$ and for all Hermitian positive definite matrices
$(a_{rs})_{r,s=1}^n$ (and all $n$) the matrix
$(f(a_{rs}))_{r,s=1}^n$ is Hermitian positive definite, then $f$
is necessarily of the form $$f(z) = \sum_{k,\ell =0}^\infty
b_{k,\ell} z^k \oz^\ell\eqno(1.6)$$ where the $b_{k,\ell}\ge 0$ for all
$k,\ell \in\ZZ_+$, and the series converges for all $z \in \CC$.
If $H$ is a complex inner product space of dimension $m$, then the
set of matrices $(<\bfz^r, \bfz^s>)_{r,s=1}^n$, all $n$, obtained
by choosing arbitrary $\bfz^r$ in $H$ is exactly the set of all
positive definite matrices of rank at most $m$.

In Section 2 we prove that this characterization (1.6) remains
valid for every complex inner product space $H$ of dimension at
least 3. That is, if dim $H\ge 3$ and the matrix
$$\left(f(<\bfz^r,\bfz^s>)\right)_{r,s=1}^n$$ is Hermitian
positive definite for all choices of $\bfz^1\nek \bfz^n$ in $H$,
all $n\in \NN$, then $f$ is necessarily of the form (1.6). For
$m=1$ this result is not valid. We do not know what happens when
$m=2$. Our proof of this result uses a generalization of the method
of proof in FitzGerald, Micchelli, Pinkus [11] and also uses the result
of Lu, Sun [19] and an extension theorem for separately real analytic
functions.

In Section 3 we assume $f$ is of the form (1.6) and
prove the following appealing characterization of
strictly Hermitian positive definite functions. Let
$$J=\{(k,\ell)\,:\, b_{k,\ell}>0\}\,.$$
We show that $f(<\cdot\,,\cdot>)$ is a strictly Hermitian positive
definite function if and only if  $(0,0)\in J$, and
$$\left| \{k-\ell:\, (k,\ell)\in J\,,\, k-\ell =q\, (\mod \,p)\}\right| =
\infty$$
for all choices of $p\in \NN$ and $q\in \{0,1\nek p-1\}$. The latter
condition simply says that
every arithmetic sequence in $\ZZ$ intersects the values
$\{k-\ell\,:\, (k,\ell)\in J\}$ an infinite number of times. Our
proof of this result utilizes the Skolem-Mahler-Lech Theorem
from number theory and a generalization
thereof due to M.~Laurent.

The paper Sun, Menegatto [30] contains some sufficient conditions for when
a function of the form (1.6) generates a strictly Hermitian positive
definite function of a fixed order $n$ on the complex unit sphere
$S^\infty$ of $\ell^2$.

\bigskip\noindent
\centerline{\bf \S 2. Hermitian Positive Definite Functions}

\medskip\noindent
In this section we prove the following.

\medskip\noindent
{\bf Theorem 2.1.} {\sl Let $H$ be a complex inner product space with
$\dim \,H\ge 3$. Then for every choice of $\bfz^1\nek
\bfz^n\in H$ the matrix
$$\left(f(<\bfz^r, \bfz^s>)\right)_{r,s=1}^n$$
is Hermitian positive definite if and only if $f$ is of the form
$$f(z) = \sum_{k,\ell =0}^\infty b_{k,\ell} z^k \oz^\ell$$
where $b_{k,\ell}\ge 0$, all $k,\ell \in\ZZ_+$, and the series
converges for all $z \in \CC$.}

\medskip
Let $\calH_m(\CC)$ denote the set of all  Hermitian positive definite
matrices of rank at most $m$, and  $\calH_m(\RR)$ the set of all
(real) positive definite matrices of rank at most $m$. If $H$ is a
complex inner product space of dimension $m$, then for any $\bfz^1\nek
\bfz^n\in H$ the matrix
$$\left(<\bfz^r, \bfz^s>\right)_{r,s=1}^n$$
is in $\calH_m(\CC)$. The converse direction also holds. That is, for
every $n$ and any  Hermitian positive definite
matrix $A=(a_{rs})_{r,s=1}^n$ of rank at most $m$, there exist
$\bfz^1\nek \bfz^n\in H$ such that
$$a_{rs} = <\bfz^r, \bfz^s>\,,\qquad r,s=1\nek n\,.$$
Thus characterizing all Hermitian positive definite functions on $H$
is equivalent to characterizing these $f:\CC\to\CC$ for which
$$\left(f(a_{rs})\right)_{r,s=1}^n$$
is Hermitian positive definite for all possible Hermitian positive definite
matrices $A=(a_{rs})_{r,s=1}^n$ of rank at most $\dim\,H$. This same
result holds over the reals.
We will denote by $\calF_m(\CC)$ the class of all Hermitian
positive definite functions on Hermitian positive definite matrices of
rank at most $m$, and by $\calF_m(\RR)$ the analogous set with respect
to the reals.

For $m=1$ the claim of Theorem 2.1 is not valid. For example, it is
readily verified that
$$f(z) =\cases{|z|^\alp,&  $z\ne 0$\cr 0, & $z=0$\cr}$$
is in $\calF_1(\CC)$ for every $\alp\in \RR$. It is possible that
Theorem 2.1 is valid for $m=2$. This case remains open.

It will be convenient to divide the proof of Theorem 2.1 into a
series of steps. Sufficiency is simple and known. We include a
proof thereof for completeness. This and the next lemma are also
to be found in FitzGerald, Micchelli, Pinkus [11].

\medskip\noindent
{\bf Lemma 2.2.} {\sl Assume $f,g\in \calF_m(\CC)$ and $a,b \ge
0$. Then
\item{(i)} $af+bg\in \calF_m(\CC)$
\item{(ii)} $f\cdot g\in \calF_m(\CC)$
\item{(iii)} the three functions $1, z, \oz$ are in $\calF_m(\CC)$
\item{(iv)} $\overline{f} \in \calF_m(\CC)$
\item{(v)} $\calF_m(\CC)$ is closed under pointwise convergence.

}

\medskip\noindent
{\bf Proof.} (i) follows from the fact that $\calH_\infty(\CC)$,
the set of all Hermitian positive definite matrices, is a positive
cone. (v) is a consequence of the closure of  $\calH_\infty(\CC)$.
If the $n\times n$ matrix $A$ is in $\calH_m(\CC)$, then the
$n\times n$ matrices $f(A)$ and $g(A)$ ($f(A) =
(f(a_{rs}))_{r,s=1}^n$) are in $\calH_\infty(\CC)$. (ii) is a
consequence of the Schur Product Theorem, see e.g. Horn, Johnson
[15, p.~309]. That is, for $A\in \calH_m(\CC)$ the matrix
$f(A)g(A) = (f(a_{rs})g(a_{rs}))_{r,s=1}^n$ is in
$\calH_\infty(\CC)$. (iii) and (iv) essentially follows by
definition. \eop

\smallskip
From Lemma 2.2 immediately follows the sufficiency part of Theorem
2.1. We now consider the necessity.

\medskip\noindent
{\bf Lemma 2.3.} {\sl Let $f\in \calF_m(\CC)$, and set
$$f(z) = u(z) + iv(z)$$
where $u(z) = {\rm Re}\, f(z)$ and $v(z) ={\rm Im}\, f(z)$. Then for $z\in \CC$
\item{(i)} $f(\oz) = \overline{f(z)}$
\item{(ii)} $u(\oz) = u(z)$
\item{(iii)} $v(\oz)=-v(z)$
\item{(iv)} $u\in \calF_m(\CC)$
\item{(v)} $f|_\RR = u|_\RR \in \calF_m(\RR)$

}

\medskip\noindent
{\bf Proof.} (i) If $A\in \calH_m(\CC)$, then $\oa_{rs} = a_{sr}$. As
$f(A)\in \calH_\infty(\CC)$ we also have
$f(a_{sr})=\overline{f(a_{rs})}$. Thus
$$f(\oa_{rs}) = f(a_{sr})=\overline{f(a_{rs})}$$
and $f(\oz) = \overline{f(z)}$ for all $z\in \CC$. (ii) and (iii) are
consequences of (i). To prove (iv) we note that from Lemma 2.2 (iv) we
have $\overline{f} \in \calF_m(\CC)$. Thus from Lemma 2.2 (i)
$$u={{f+\overline{f}}\over 2}\in \calF_m(\CC)\,.$$
To prove (v) note that $f|_\RR = u|_\RR$ as a consequence of
(iii). Now $u|_\RR:\RR\to\RR$ and $u\in \calF_m(\CC)$. Thus
$ u|_\RR \in \calF_m(\RR)$. \eop

\medskip\noindent
{\bf Lemma 2.4.} {\sl Let $f\in \calF_m(\CC)$, $m\ge 2$.
Then for any $a, b\in \CC$
$$g(z) =f(z+1) + f(|a|^2z+|b|^2) \pm \left[ f(a z+b) + f(\oa z+
\ob)\right]$$
is in $\calF_{m-1}(\CC)$.}

\medskip\noindent
{\bf Proof.} Assume $A$ is an $n\times n$ matrix in
$\calH_{m-1}(\CC)$. Then for any $a\in \CC$ the $2n\times 2n $
matrix
$$\left(\matrix{ A & a A\cr  \oa A & |a|^2A\cr}\right)$$
is also in $\calH_{m-1}(\CC)$.

Let $J$ denote the $n\times n$ matrix all of whose entries
are one. $J\in \calH_1(\CC)$, and the $2n \times 2n$ matrix
$$\left( \matrix{ J & b J\cr
\ob J & |b|^2 J\cr}\right)$$
is also in $\calH_1(\CC)$ for any $b\in \CC$. Thus
$$\left(\matrix{ A + J  & aA + bJ \cr
\oa A +\ob J  & |a|^2 A + |b|^2 J \cr}\right)$$
is in $\calH_m(\CC)$. As $f\in \calF_m(\CC)$ it follows that
$$\left(\matrix{ f(A + J)  & f(aA +bJ)
\cr  f(\oa A +\ob J)  & f(|a|^2 A + |b|^2 J) \cr}\right)$$
is in $\calH_\infty(\CC)$.

For any $\bfc^T=(c_1\nek c_n)\in \CC^n$, set
$$(\bfc^T, \pm\bfc^T) =(c_1\nek c_n, \pm c_1\nek \pm c_n)\in
\CC^{2n}\,.$$
Thus
$$(\bfc^T, \pm\bfc^T)
\left(\matrix{ f(A + J)  & f(a A + b J)
\cr  f(\oa A +\ob J)  & f(|a|^2 A + |b|^2 J) \cr}\right)
\left(\matrix{\bfc \cr \pm\bfc\cr}\right) \ge 0\,.$$
That is,
$$\bfc^T\left[  f(A + J)  +  f(|a|^2 A + |b|^2 J) \pm
\left[ f(a A + b J) +  f(\oa A
+\ob J) \right]\right]\bfc\ge 0$$
which implies that
$$g(z) =f(z+1) + f(|a|^2 z+|b|^2) \pm \left[ f(a z+ b)
+  f(\oa z+ \ob)\right]$$
is in $\calF_{m-1}(\CC)$. \eop

Setting $z=t\in \RR$ it follows from Lemma 2.4 and Lemma  2.3 (v)
that
$$h(t) =u(t+ 1) + u(|a|^2 t+|b|^2) \pm \left[ u(a t+b)
+  u(\oa t+ \ob)\right]$$
is in $\calF_{m-1}(\RR)$. Furthermore from Lemma 2.3 (ii)
$$u(\oa t + \ob) = u(a t + b)\,.$$
Thus
$$h(t) = u(t+ 1) + u(|a|^2 t+|b|^2) \pm 2 u(a t+b)$$
is in $\calF_{m-1}(\RR)$ for any $a, b\in \CC$.

We shall use the following result.

\medskip\noindent
{\bf Theorem 2.5.} (Lu, Sun [19]) {\sl Assume $f\in \calF_m(\RR)$,
$m\ge 2$. Then $$f(t) =\sum_{k=0}^\infty c_k t^k$$ with $c_k\ge 0$
for all $k$, and where the series converges for all $t\in \RR$.}

\medskip
There is a slight oversight in the proof of Lu, Sun [19]. It is
also necessary in their proof that $f$ be bounded in a
neighborhood of the origin. It may be easily shown that this holds
for $f\in \calF_m(\RR)$ when $m\ge 2$.

Thus for $m\ge 3$ we have that for all $a, b\in \CC$
$$h(t) =u(t+1) + u(|a|^2 t+|b|^2) \pm 2 u(a t+ b)\eqno(2.1)$$
is in $\calF_2(\RR)$ and thus has a power series expansion, with
nonnegative coefficients, which converges for all $t$. As
$u(t)$ has this property, so does
$u(t+1)$ and $u(|a|^2 t+|b|^2)$. Thus
$$u(at+b)$$
has a power series expansion in $t$, which converges for all $t\in \RR$,
for each fixed $a,b\in\CC$. In this power series expansion the
coefficients need not be nonnegative.

It will be convenient to consider $u$ as a map from $\RR^2$ to
$\RR$, rather than from $\CC$ to $\RR$. Thus if $w=w_1+iw_2 \in \CC$, we
set $\bfw=(w_1,w_2)\in \RR^2$ and write
$$U(\bfw) = u(w)\,.$$
However we also write
$$U(\bfa t +\bfb)= u(at+b)$$
as we will consider $t\in \RR$ as a parameter.

\medskip\noindent
{\bf Proposition 2.6.} {\sl Let $U$ be as above and $m\ge 3$.
Then
\item{(i)} For each $\bfa, \bfb \in \RR^2$ we have
$$U(\bfa t+\bfb) = \sum_{k=0}^\infty c_k(\bfa, \bfb)t^k
\eqno(2.2)$$
where the power series converges for all $t\in \RR$.
\item{(ii)} $U\in C(\RR^2)$.
\item{(iii)} For each $M> 0$ there exists a sequence of
positive numbers $(b_k(M))_{k=0}^\infty$ such that
for each $k$, and all $\bfa, \bfb\in \RR^2$ satisfying
$|\bfa|, |\bfb|\le M$ we have
$$|c_k(\bfa, \bfb)|\le b_k(M)$$
\item{}and}
$$\limsup_{k\to\infty} \root k \of{b_k(M)}=0\,.$$

\medskip\noindent
{\bf Proof.} (i) is just a restatement of the consequence
of Theorem 2.5. We prove (ii) as follows. Let $h$ be as in (2.1).
Since $h\in \calF_2(\RR)$, it follows that
$$h(t) =\sum_{k=0}^\infty c_k t^k$$
where $c_k \ge 0$ and the series converges for all $t \in\RR$. Thus
$h$ is also increasing and nonnegative on $\RR_+$. We therefore
have $h(t) -h(0)\ge 0$
for all $t\in \RR_+$ that can be rewritten as
$$u(t+1) - u(1) + u(|a|^2 t+|b|^2) - u(|b|^2)  \ge
2|u(at+b) - u(b)|\,.$$
Fix $b\in \CC$. As $u|_\RR\in C(\RR)$, it follows that given $\eps>0$,
we have a $\delta>0$ such that for every $t\in \RR$ satisfying
$|t|<\delta$ and $a\in \CC$, $|a|\le 1$,
$$|u(t+1) - u(1)| <\eps$$
and
$$ |u(|a|^2 t+|b|^2) - u(|b|^2)|\le |u(t+|b|^2)-u(|b|^2)|
<\eps\,.$$
Thus,
$$|u(at+b)-u(b)|<\eps$$
for all $a\in \CC$ satisfying $|a|\le 1$ and $t\in [0, \delta)$.
This proves (ii).

From the above
$$u(t+ 1) + u(|a|^2 t+|b|^2) \pm 2 u(a t+b)$$
has a power expansion with nonnegative coefficients, as does
$u(t+ 1)$ and $ u(|a|^2 t+|b|^2)$. Furthermore if
$$u(|a|^2 t+|b|^2) = \sum_{k=0}^\infty d_k(|a|^2, |b|^2) t^k$$
and $|a|, |b| \le M$, then as is easily verified
$$0\le d_k(|a|^2, |b|^2) \le d_k(M^2, M^2)\,.$$
Thus
$$|c_k(\bfa, \bfb)| \le d_k(1,1) + d_k(M^2, M^2)\,.$$
for each $k$. As the
$$\sum_{k=0}^\infty d_k(M^2, M^2) t^k$$
is entire we have
$$\limsup_{k\to\infty} \root k \of{d_k(M^2, M^2)} =0$$
for each $M$. This proves (iii). \eop

\medskip
As a consequence of Proposition 2.6 we have the following result.

\medskip\noindent
{\bf Theorem 2.7.} {\sl Assume $U$ satisfies the conditions of
Proposition 2.6. Then $U$ is the
restriction to $\RR^2$ of an entire function $\tilU$ on $\CC^2$.}

\medskip
To verify Theorem 2.7 we need less than we have proven in Proposition
2.6. In (i) of Proposition 2.6 we proved that $U(\bfw)$ has a
convergent power series expansion on every straight line in $\RR^2$.
It suffices, in the proof of Theorem 2.7, for this property to only hold on
lines parallel to the axes. From that and the other results proved
in Proposition 2.6 follows Theorem 2.7 as a consequence of a
result in Bernstein [3, p.~101]. It also follows from results
in Browder [5], Cameron, Storvick [6] and Siciak [29], see
the review
article by Nguyen [22]. Theorem 2.7 implies that the series
expansion for $U$ converges appropriately. That is, $$U((x,y)) =
u(x+iy) = \sum_{k,\ell=0}^\infty a_{k,\ell} x^k y^\ell$$ where the
power series converges absolutely for all $(x,y)\in \RR^2$.

We now continue as in FitzGerald, Micchelli, Pinkus [11].

\medskip\noindent
{\bf Proposition 2.8.} {\sl Let $f\in \calF_m(\CC)$, $m\ge 3$, and
$v=\,{\rm Im\,}f$. Then
$$v(x+iy) = \sum_{k,\ell=0}^\infty c_{k,\ell} x^k y^\ell$$
where the power series converges absolutely for all $(x,y)\in \RR^2$.}

\medskip\noindent
{\bf Proof.} From Lemma 2.2 (ii) and (iii) we have $zf(z)\in
\calF_m(\CC)$. Thus from Theorem  2.7 its real part
$$xu(x+iy) - yv(x+iy)$$
has a power series expansion, which in turn implies that
$$yv(x+iy) =  \sum_{k,\ell=0}^\infty d_{k,\ell} x^k y^\ell$$
with the power series converging absolutely. Setting $y=0$ we obtain
$$0 = \sum_{k,\ell=0}^\infty d_{k,0} x^k\,.$$
Thus $d_{k,0}=0$ for all $k$ and
$$v(x+iy) =  \sum_{k,\ell=0}^\infty d_{k,\ell+1} x^k y^\ell\,,$$
where the power series converges absolutely. \eop

\smallskip
We can now finally prove Theorem 2.1.

\medskip\noindent
{\bf Proof of Theorem 2.1.} As $f\in \calF_m(\CC)$, $m\ge 3$,
we have from Theorem 2.7 and Propositions 2.8 that
$$f(x+iy) =  \sum_{k,\ell=0}^\infty e_{k,\ell} x^k y^\ell$$
and the power series converges absolutely for all $(x,y)\in
\RR^2$. Substituting $x=(z+\oz)/2$ and $y=(z-\oz)/2$ we obtain
$$f(z) =  \sum_{k,\ell=0}^\infty b_{k,\ell} z^k \oz^\ell$$
and this  power series also converges absolutely for all $z\in
\CC$. It remains to prove that $b_{k,\ell} \ge 0$
for all $(k,\ell)\in
\ZZ_+^2$. For $\eps>0$ and $z_1\nek z_n\in \CC$, the matrix
$$A=(\eps z_r\oz_s)_{r,s=1}^n$$
is Hermitian positive definite of rank 1. Thus
$$\sum_{r,s=1}^n c_r f(\eps z_r\oz_s)\oc_s\ge 0$$
for all $c_1\nek c_n \in \CC$. Substituting for $f$ we obtain
$$\sum_{k,\ell=0}^\infty b_{k,\ell} \eps^{k+\ell} \sum_{r,s=1}^n c_r
z_r^k \oz_s^k \oz_r^\ell z_s^\ell \oc_s = \sum_{k,\ell=0}^\infty
b_{k,\ell} \eps^{k+\ell} \left| \sum_{r=1}^n c_r z_r^k \oz_r^\ell
\right|^2\ge 0\,.$$
This inequality must hold for all choices of $n\in \NN$, $\eps>0$, and
$z_1\nek z_n, c_1\nek c_n \in \CC$. Given $(i,j)\in \ZZ_+^2$ it is
possible to choose $n$ and $z_1\nek z_n, c_1\nek c_n \in \CC$ such
that
$$\sum_{r=1}^n c_r z_r^k \oz_r^\ell=0$$
for all $(k,\ell)\in \ZZ_+^2$ satisfying $k+\ell \le i+j$
except that
$$\sum_{r=1}^n c_r z_r^i \oz_r^j \ne 0\,.$$
Letting $\eps \downarrow 0$ then proves $b_{i,j}\ge 0$. \eop

\bigskip\noindent
\centerline{\bf \S 3. Strictly Hermitian Positive Definite Functions}

\medskip\noindent
In this section we always assume that $f$ is of the form
$$f(z) = \sum_{k,\ell =0}^\infty b_{k,\ell} z^k \oz^\ell\eqno(3.1)$$
where the $b_{k,\ell}\ge 0$ for all $k,\ell$ in $\ZZ_+$, and the series
converges for all $z \in \CC$. We prove the following result.

\medskip\noindent
{\bf Theorem 3.1.} {\sl Assume $f$ is of the form (3.1). Set
$$J=\{(k,\ell)\,:\, b_{k,\ell}>0\}\,.$$
Then $f$ is strictly Hermitian positive definite on $H$ if and only if
$(0,0)\in J$, and for each $p\in \NN$ and $q\in \{0,1\nek p-1\}$}
$$\left| \{k-\ell:\, (k,\ell)\in J\,,\, k-\ell =q\, (\mod \,p)\}\right| =
\infty\,.$$

\medskip
The conditions of Theorem 3.1 are independent of $H$. We first
show that it suffices in the proof of Theorem 3.1 to only consider
the standard inner product on $\CC$. In this we follow the
analysis in Pinkus [23]

\medskip\noindent
{\bf Proposition 3.2.} {\sl Theorem 3.1 is valid if and only if it
holds for the standard inner product on $\CC$, namely
$$<z, w>= z\ow\,,$$
for $z,w\in \CC$.}

\medskip
The main tool used in the proof of Proposition 3.2 is Proposition
3.3 which appears in Pinkus [23] in the real case, and is
essentially based on an exercise in P\'olya, Szeg\"o [24, p.~287].

\medskip\noindent
{\bf Proposition 3.3.} {\sl  Let $H$ be a complex inner product space, and
$\bfz^1\nek \bfz^n$ any $n$ distinct points in $H$. There then exist
distinct $z_1\nek z_n\in \CC$ and a Hermitian positive definite matrix
$(m_{rs})_{r,s=1}^n$ such that}
$$<\bfz^r,\bfz^s> = z_r \oz_s + m_{rs}\,.$$

\medskip\noindent
{\bf Proof.} Set
$$a_{rs} = <\bfz^r,\bfz^s>\,,\qquad r,s=1\nek n\,.$$
Since the $\bfz^1\nek \bfz^n$ are distinct points in
$H$, the Hermitian positive definite matrix
$$A=(a_{rs})_{r,s=1}^n$$
has no two identical rows (or columns). For assume there are two
identical rows indexed by $i$ and $j$. Then
$$a_{ii} = a_{ij} = a_{ji} =a_{jj}$$
implying that
$$ <\bfz^i,\bfz^i>=<\bfz^i,\bfz^j>= <\bfz^j,\bfz^j>\,.$$
But then
$$\|  <\bfz^j,\bfz^j> \bfz^i -  <\bfz^i,\bfz^i> \bfz^j\|=0$$
and thus $\bfz^i = \bfz^j$, contradicting our assumption.

As $A$ is an $n\times n$ Hermitian positive definite matrix it may
be decomposed as
$$A= C^T \oC$$
where $C=(c_{kr})_{k=1}^m{}_{r=1}^n$ is an $m\times n$ matrix,
$\oC=(\oc_{kr})_{k=1}^m{}_{r=1}^n$
and $\rank A=m$. Thus
$$a_{rs} = \sum_{k=1}^m c_{kr} \oc_{ks}\,,\qquad r,s=1\nek n\,.$$
Since no two rows of $A$ are identical, it follows that no two
columns of $C$ are identical. As such there exists a $\bfv\in
\CC^m$, $\|\bfv\|=1$, for which
$$\bfv C = \bfz$$
with
$\bfz=(z_1\nek z_n)$ where the $z_1\nek z_n$ are all distinct.

Let $V$ be any $m\times m$ unitary matrix ($V^T\oV=I$) whose
first row is the above $\bfv$. Let $U$ be the $m\times m$ matrix
whose first row is the above $\bfv$ and all its other entries are
zero. Then
$$A= C^T\oC = C^TV^T\oV\oC= C^TU^T\oU\oC +C^T(V-U)^T(\oV-\oU)\oC$$
since $U^T(\oV-\oU) = (V-U)^T\oU = 0$. From the above it follows that
$(C^TU^T\oU\oC)_{rs} = z_r\oz_s$ and $M= C^T(V-U)^T(\oV-\oU)\oC$ is
Hermitian positive definite. This proves the proposition. \eop

\smallskip
We now prove Proposition 3.2.

\medskip\noindent
{\bf Proof of Proposition 3.2.} Let $H$ be any complex inner product
space, and let $\bfz\in H$, $\|\bfz\|=1$. Then for $z, w\in
\CC$
$$f(<z\bfz, w\bfz>) = f(z\ow)\,.$$
This immediately implies that if $f$ is not a strictly Hermitian
positive definite function on $\CC$, then it is not a strictly Hermitian
positive definite function on any complex inner product space $H$.

The converse direction is a consequence of Proposition 3.3. Assume
$f$ is a strictly Hermitian positive definite function on $\CC$.
Let $\bfz^1\nek \bfz^n$ be distinct points in $H$. By Proposition
3.2 there exist distinct $z_1\nek z_n\in \CC$ and a Hermitian
positive definite matrix $(m_{rs})_{r,s=1}^n$ such that
$$<\bfz^r,\bfz^s> = z_r \oz_s + m_{rs}\,.$$ As the
$(m_{rs})_{r,s=1}^n$ is a Hermitian positive definite matrix, we
have for any nonzero $c_1\nek c_n$ $$\sum_{r,s=1}^n c_r
<\bfz^r,\bfz^s>\oc_s = \sum_{r,s=1}^n c_r z_r \oz_s\oc_s +
\sum_{r,s=1}^n c_r m_{rs} \oc_s\ge \sum_{r,s=1}^n c_r z_r
\oz_s\oc_s \,.$$ Similarly, it is well known that $$\sum_{r,s=1}^n
c_r <\bfz^r,\bfz^s>^k \overline{<\bfz^r,\bfz^s>}^\ell \oc_s \ge
\sum_{r,s=1}^n c_r (z_r \oz_s)^k  (\oz_r z_s)^\ell \oc_s\,.$$ This
is a consequence of inequalities between the Schur (Hadamard)
product of two Hermitian positive definite matrices, see e.g.,
Horn, Johnson [15, p.~310]. Thus for $f$ of the form (3.1)
$$\sum_{r,s=1}^n c_r f(<\bfz^r,\bfz^s>) \oc_s = \sum_{k,
\ell=0}^\infty b_{k,\ell} \sum_{r,s=1}^n c_r <\bfz^r,\bfz^s>^k
<\bfz^s,\bfz^r>^\ell \oc_j$$
$$\ge \sum_{k,\ell=0}^\infty b_{k, \ell} \sum_{r,s=1}^n c_r (z_r\oz_s)^k
(\oz_r z_s)^\ell\oc_s = \sum_{r,s=1}^n c_r f(z_r \oz_s) \oc_s
\,.$$
As the $z_1\nek z_n$ are distinct, the $c_1\nek c_n$ are nonzero
and $f$ is a  strictly Hermitian positive definite function on $\CC$,
it follows that the above quantity is strictly positive. Thus
$f$ is a strictly Hermitian positive definite function on $H$.
\eop

\medskip
We therefore assume in what follows that $H=\CC$. Let $z_1\nek
z_n$ be $n$ distinct points in $\CC$ and $c_1\nek c_n \in \CC\\
\{0\}$. We want conditions implying
$$\sum_{r,s=1}^n c_r f(z_r\oz_s)\oc_s>0\,.$$
From the form of $f$, i.e., (3.1), we have
$$\sum_{r,s=1}^n c_r f(z_r\oz_s)\oc_s = \sum_{r,s=1}^n c_r
\left(\sum_{k,\ell=0}^\infty b_{k,\ell}(z_r\oz_s)^k (\oz_r
  z_s)^\ell\right)
\oc_s =
\sum_{k,\ell=0}^\infty b_{k,\ell}\left|\sum_{r=1}^n c_rz_r^k
  \oz_r^\ell \right|^2\,.$$
Now $b_{k,\ell}\ge 0$ and $|\sum_{r=1}^n c_rz_r^k \oz_r^\ell|^2\ge
0$. Thus $f$ is a strictly Hermitian positive definite function if and
only if for all choices of $n$, distinct points $z_1\nek z_n$ in $\CC$
and nonzero values $c_1\nek c_n$ we always have
$$\sum_{r=1}^n c_rz_r^k \oz_r^\ell \ne 0\eqno(3.2)$$
for some $(k,\ell)\in J$.

One direction in the proof of Theorem 3.1 is elementary.

\medskip\noindent
{\bf Proposition 3.4.} {\sl If $f$  is a strictly Hermitian positive
definite function then $(0,0)\in J$ and for each $p\in \NN$, $q\in
\{0,1\nek p-1\}$ we have}
$$\left|\{ k-\ell\,:\,(k,\ell)\in J, k-\ell = q\,(\mod\,p)\}\right|
=\infty\,.$$

\smallskip\noindent
{\bf Proof.} To prove that we must have $(0,0)\in J$ we simply take
$n=1$, $z_1=0$ and $c_1\ne 0$ in (3.2). As $f$ is  a strictly
Hermitian positive  definite function then
$$c_1z_1^k\oz_1^\ell \ne 0$$
for some $(k,\ell)\in J$. But
$$c_1z_1^k\oz_1^\ell = 0$$
for all $(k,\ell)\in \ZZ_+^2\\ \{(0,0)\}$, which implies that we must
have $(0,0)\in J$.

Let us now assume that there exists a $p\in \NN$ and $q\in \{0,1\nek
p-1\}$ for which
$$\left|\{ k-\ell\,:\,(k,\ell)\in J, k-\ell = q\,(\mod\,p)\}\right| =N
<\infty\,.$$
Thus if $(k,\ell)\in J$ and  $k-\ell = q\,(\mod\,p)$, then
$$k-\ell = q + a_m p\,,\qquad m=1\nek N\,,$$
for some $a_1\nek a_N \in \ZZ$.

Choose $\tet_1\nek \tet_{N+1}$ in $[0,2\pi)$ such that the $p(N+1)$
points
$$z_{r,t} = e^{i(\tet_r + {{2\pi t}\over p})}\,,\qquad r=1\nek N+1;
\,  t=0,1\nek p-1\,,$$
are all distinct.

There exist $c_1\nek c_{N+1}$, not all zero, such that
$$\sum_{r=1}^{N+1} c_r e^{i\tet_r(q+a_mp)}=0\,,\qquad m=1\nek N\,.\eqno(3.3)$$
Since the $e^{i{{2\pi t}\over p}}$, $t=0,1\nek p-1$, are distinct, the
$p\times p$ matrix $\left(e^{i{{2\pi ts}\over
p}}\right)_{t,s=0}^{p-1}$ is nonsingular and there exist
$d_0\nek d_{p-1}$, not all zero, such that
$$\sum_{t=0}^{p-1} d_t e^{i{{2\pi ts}\over p}}=0\,, \qquad s=
0,1\nek p-1\,;\, s\ne q\,.\eqno(3.4)$$
We claim that
$$\sum_{t=0}^{p-1} \sum_{r=1}^{N+1} d_t c_r z_{r,t}^k
\oz_{r,t}^\ell=0$$
for all $(k,\ell)\in J$. As the $z_{r,t}$ are distinct and the
$d_tc_r$ are not all zero, this would imply that $f$ is not a
strictly Hermitian positive definite function.

Now
$$\sum_{t=0}^{p-1} \sum_{r=1}^{N+1} d_t c_r z_{r,t}^k
\oz_{r,t}^\ell = \sum_{t=0}^{p-1} \sum_{r=1}^{N+1} d_t c_r
e^{i(\tet_r + {{2\pi t}\over p})k} e^{-i(\tet_r + {{2\pi t}\over
p})\ell}$$
$$ = \left( \sum_{t=0}^{p-1} d_t e^{i{{2\pi t}\over p}(k-\ell)}\right)
\left( \sum_{r=1}^{N+1} c_r
e^{i\tet_r(k-\ell)}\right)\,.$$
If $(k,\ell)\in J$ and $k-\ell =q\,(\mod \,p)$, then $k-\ell = q+a_mp$
and from (3.3) the right-hand factor is zero. If $(k,\ell)\in J$ and
$k-\ell \ne q\,(\mod \,p)$, then $k-\ell = s\,(\mod\, p)$, $s\in
\{0,1\nek p-1\}\\ \{q\}$. Thus
$$k-\ell = s+ mp$$
for some  $s\in \{0,1\nek p-1\}\\ \{q\}$ and $m\in \ZZ$, and
$$e^{i{{2\pi t}\over p}(k-\ell)}= e^{i{{2\pi ts}\over p}}e^{i2\pi tm}
=  e^{i{{2\pi ts}\over p}}\,.$$
From (3.4) the left-hand factor is zero. This proves the
proposition. \eop

\medskip
It is the converse direction which is less elementary. We will prove
the following result which completes the proof of Theorem 3.1.

\medskip\noindent
{\bf Theorem 3.5.} {\sl Let $J\subseteq \ZZ_+^2$ be such
that for each $p\in \NN$ and $q\in \{0,1\nek p-1\}$ we have
$$\left| \{k-\ell\,:\, (k,\ell)\in J,\, k-\ell =q\,(\mod\, p)\}\right|
= \infty\,.$$
Then for all $n$, distinct nonzero points $z_1\nek z_n$ in $\CC$ and
nonzero values $c_1\nek c_n$ in $\CC$, we always have
$$\sum_{r=1}^n c_r z_r^k\oz_r^\ell \ne 0$$
for some $(k,\ell)\in J$.}

\medskip
Note that we have here assumed that each of the $z_1\nek z_n$ is
nonzero and we have dropped the condition $(0,0)\in J$.
It is easily shown that we can make this assumption.

A {\sl linear recurrence relation} is a series of equations of the
form
$$a_{s+m}= a_{s+m-1}w_1+\cdots + a_s w_m\,,\qquad s\in \ZZ\eqno(3.5)$$
satisfied by the {\sl recurrence sequence} $\{a_s\}_{s\in \ZZ}$ for
some given $w_1\nek w_m$ ($m$ finite). We assume
that the $w_j$ and $a_s$ are in $\CC$. Associated with each such
recurrence sequence is a {\sl generalized power sum}
$$a_s =\sum_{j=1}^r P_j(s) u_j^s\eqno(3.6)$$
where the $P_j$ are polynomials in $s$ of degree $\partial P_j$, with
$\sum_{j=1}^r \partial P_j +1 \le m$, and vice versa. That is, each
generalized power sum of the form (3.6) gives rise to a linear
recurrence relation of the form (3.5).

We are interested in the form (3.6). A delightful theorem in
number theory which concerns recurrence sequences is the
Skolem-Mahler-Lech Theorem, see e.g. Shorey, Tijdeman [28, p.~38].

\medskip\noindent
{\bf Theorem 3.6.} (Skolem--Mahler-Lech) {\sl Assume that
$\{a_s\}_{s\in\ZZ}$ is a recurrence sequence, i.e., satisfies (3.5) or
(3.6). Set
$$\calA = \{s\,:\,a_s=0\}\,.$$
Then $\calA$ is the union of a finite number of points and a finite
number of full arithmetic sequences.}

\medskip
What is the connection between this result and our problem? Recall
that we are concerned with the equations
$$\sum_{r=1}^n c_r z_r^k\oz_r^\ell \,.$$
Assume for the moment that $|z_r|=1$ for all $r=1\nek n$. We can
then rewrite the above as
$$\sum_{r=1}^n c_r z_r^{k-\ell}\,.$$
Set
$$a_s = \sum_{r=1}^n c_r z_r^s\,,\qquad s\in\ZZ\,.$$
These are equations of the form (3.6) and thus
the $\{a_s\}_{s\in \ZZ}$ is a recurrence sequence. Let us rewrite and
prove Theorem 3.5 in this particular case. (The same holds if we
assume $|z_r|=\lam$ for some $\lam\in\RR_+$ and all $r=1\nek n$.)

\medskip\noindent
{\bf Proposition 3.7.} {\sl Let $J^*\subseteq \ZZ$ be such that
for each $p\in \NN$ and $q\in \{0,1\nek p-1\}$ we have
$$\left| \{s\,:\, s\in J^*,\, s =q\,(\mod\, p)\}\right|
= \infty\,.$$
Then for all $n$, distinct nonzero points $z_1\nek z_n$ in $\CC$ and
nonzero values $c_1\nek c_n$ in $\CC$, we always have
$$\sum_{r=1}^n c_r z_r^s \ne 0$$
for some $s\in J^*$.}

\medskip\noindent
{\bf Proof.} Set
$$a_s = \sum_{r=1}^n c_r z_r^s\,,\qquad s\in\ZZ\,,$$
and let
$$\calA = \{s\,:\,a_s=0\}\,.$$
Assume $a_s=0$ for all $s\in J^*$, i.e., $J^*\subseteq \calA$. A
contradiction immediately ensues from the Skolem-Mahler-Lech Theorem
3.6 since it implies that $J^*$ is contained in the set $\calA$
which is the union of a finite number of points and a finite
number of full arithmetic sequences. As the $a_s$ cannot
all be zero, it is readily verified that there must exist an
arithmetic sequence disjoint from
$\calA$. Thus there exists a $p\in \NN$ and $q\in \{0,1\nek p-1\}$ for which
$$\{s\,:\, s\in J^*,\, s =q\,(\mod\, p)\}= \emptyset\,.$$
This proves the proposition. \eop

\noindent
{\bf Remark 1.} The above is equivalent to the following. Let
$J^*\subseteq \ZZ$ and
$$\Pi_{J^*} =\span\{ z^k\, :\, k\in J^*\}\,.$$
Then for every finite point set $E$ in $\CC$
$$\dim \Pi_{J^*}|_E = |E|$$
if and only if $0\in J^*$ and $J^*$ satisfies the criteria of
Proposition 3.7.

\medskip\noindent
{\bf Remark 2.} By a {\sl discrete measure} $d\mu$ on $[0,2\pi)$ we
mean a measure of the form
$$d\mu = \sum_{r=1}^n c_r \delta_{\tet_r}$$
for some finite $n$, where $\delta_\tet$ is the Dirac-Delta
point measure at $\tet$. From Propositions 3.4 and 3.7 the
$J^*\subseteq \ZZ$, as above, characterize the sets of uniqueness
for the Fourier coefficients of a discrete measure. Set
$$\widehat{\mu}(s) =  \int_0^{2\pi} e^{i\tet s}d\mu(\tet)\,,
\qquad s\in \ZZ\,.$$
If $d\mu_1$ and $d\mu_2$ are two discrete measures and
$$\widehat{\mu}_1(s) = \widehat{\mu}_2(s)$$
for all $s\in J^*$, then $d\mu_1=d\mu_2$.

\medskip
Unfortunately there seems to be no a priori reason to assume
that the $\{z_r\}$
are all of equal modulus. We will prove nonetheless that such is
effectively the case. We prove this fact via the following
generalization of the Skolem-Mahler-Lech Theorem.

\medskip\noindent
{\bf Theorem 3.8.} (Laurent [17], [18, p.~26]) {\sl Let
$$a_{k,\ell} = \sum_{r=1}^n c_r z_r^k w_r^\ell$$ for some $(z_r,
w_r)\in \CC^2$ and  $c_r\in \CC$, $r=1\nek n$. Let $$\calB = \{
(k,\ell)\,:\, a_{k,\ell}=0\}\,.$$ Then $\calB$ is the union of a
finite number of translates of subgroups of $\ZZ^2$.}

\medskip
Note that this is not a full generalization of the Skolem-Mahler-Lech
Theorem in that the polynomial parts of the generalized power sum, see
(3.6), are here assumed to be constants.

In our problem we are given
$$a_{k,\ell} = \sum_{r=1}^n c_r z_r^k \oz_r^\ell$$
i.e., $\oz_r=w_r$. Additionally, and importantly, we only
consider $(k,\ell)\in \ZZ_+^2$ although, since $z_r\ne 0$ for all
$n$, we can and do define $a_{k,\ell}$ for all $(k,\ell)\in \ZZ^2$.

What are the subgroups of $\ZZ^2$? Each subgroup is given
as the set of $(k,\ell)$ satisfying
$$\left(\matrix{k\cr \ell\cr}\right) = p\left(\matrix{a\cr b\cr}\right) +
q\left(\matrix{c\cr d\cr}\right)$$
where $p, q$ vary over all $\ZZ$, and $a,b,c,d$ are given integers
with $ad-bc\not\in \{-1, 1\}$. If $ad-bc\in \{-1, 1\}$, then the set
of solutions is all of $\ZZ^2$.
If $ad-bc=0$, then the solution set can be rewritten as
$$(k,\ell) = (js, jt)\,,\qquad j\in \ZZ$$
for some $(s,t)\in \ZZ^2$. This includes the case $(s,t) =(0,0)$. For
$|ad-bc|\ge 2$, the subgroup is a lattice.

Let us denote a translate of a subgroup of $\ZZ^2$ by $L$. We will
prove the following.

\medskip\noindent
{\bf Proposition 3.9.} {\sl Assume $|L\cap \ZZ_+^2|=\infty$, and
$$\sum_{r=1}^n c_r z_r^k \oz_r^\ell=0$$
for all $(k,\ell)\in L$. Then
$$\sum_{\{r: |z_r|=\lam\}} c_r z_r^k \oz_r^\ell=0$$
for all $(k,\ell)\in L$ and each $\lam\in \RR_+$.}

\medskip
Of course we only consider values $\lam$ equal  to one of
the $|z_1|\nek |z_n|$. For other values of $\lam$ the set
$\{r: |z_r|=\lam\}$ is empty.

Before proving Proposition 3.9 we first prove an ancillary result.

\medskip\noindent
{\bf Proposition 3.10.} {\sl Assume
$$\sum_{r=1}^n d_r w_r^j=0$$
for all $j\in \ZZ$ with $d_r\ne 0$ and $w_r\in \CC\\ \{0\}$, $r=1\nek
n$. Then
$$\sum_{\{r: w_r=\mu\}} d_r w_r^j=0$$
for each $\mu \in \CC$ and all $j\in\ZZ$.}

\medskip\noindent
{\bf Proof.} Let $\mu_1\nek \mu_m$ denote the distinct values of the
$w_1\nek w_n$, and assume
that $\ell_s$ of the $w_j$'s equal $\mu_s$, $s=1\nek m$. Thus
$\sum_{s=1}^m \ell_s=n$.

For any $k$
$$\det \left(w_r^{j+k}\right)_{r=1\, j=0}^{n\phantom{12} n-1} = w_1^k\cdots
w_n^k \prod_{1\le t<s\le n} (w_s-w_t)\,.$$
Thus
$$\rank \left(w_r^j\right)_{r=1,\,j\in \ZZ}^{n\phantom{123}} = m\,.$$
The linear subspace
$$D=\{\bfd=(d_1\nek d_n) \,:\, \sum_{r=1}^n d_r w_r^j=0\,,\ {\rm all\ }j\in \ZZ\}$$
is therefore of dimension $n-m$.

Now if, for example, $w_1=\cdots
=w_{\ell_1}=\mu_1$ (and $w_j\ne \mu_1$ for $j>\ell_1$), then
$$\rank \left( w^j_t\right)_{t=1,\,j\in \ZZ}^{\ell_1\phantom{123}}
=1$$
and there are therefore $\ell_1-1$ linearly independent
vectors $\bfd^1\nek \bfd^{\ell_1}$ in $D$ satisfying
$$d_t^1=\cdots =d_t^{\ell_1} =0$$
for all $t=\ell_1+1 \nek n$. In this manner we obtain
$$\sum_{s=1}^m \left(\ell_s-1\right) = n-m$$
linearly independent vectors in $D$. As
$$\dim D = n-m$$
we have therefore found a basis for $D$. The
proposition now easily follows. \eop

\medskip\noindent
{\bf Proof of Proposition 3.9.} A translate of a subgroup $L$ of
$\ZZ^2$ is given by the formula
$$\left(\matrix{k\cr \ell\cr}\right) = \left(\matrix{\alp\cr \beta\cr}\right)
+ p\left(\matrix{a\cr b\cr}\right) +
q\left(\matrix{c\cr d\cr}\right)$$
where $p, q$ vary over all $\ZZ$. We divide the proof of this
proposition into two cases.

\smallskip\noindent
{\bf Case 1.} $ad-bc=0$. In this case
$$(k,\ell) = (\alp + js, \beta +jt)\,.$$
If $st<0$ or if $s=t=0$, then $|L\cap \ZZ_+^2|<\infty$. As such, we may
assume that $s, t\ge 0$ and $s+t > 0$. Now for $(k,\ell)\in L$
$$0=\sum_{r=1}^n c_r z_r^k \oz_r^\ell = \sum_{r=1}^n c_r z_r^{\alp+js}
\oz_r^{\beta +jt}
= \sum_{r=1}^n \left(c_r z_r^\alp \oz_r^\beta\right)\left(z_r^s
\oz_r^t\right)^j = \sum_{r=1}^n d_r w_r^j$$
for all $j\in \ZZ$, where
$$d_r = c_r z_r^\alp \oz_r^\beta$$
and
$$ w_r = z_r^s \oz_r^t\,.$$
Note that $d_r\ne 0$ since $c_r\ne 0$ and $z_r\ne 0$.
Furthermore as $s+t>0$, and
$$|w_r| = |z_r|^{s+t}$$
we have
$$\{r\,: \,|z_r|=\lam\} =\{r\,:\, |w_r|=\lam^{s+t}\}$$
i.e., the indices $\{1\nek n\}$ divide in the same way when
considering the distinct $|z_r|$ or the distinct $|w_r|$. We apply
Proposition 3.10 to obtain our result.

\smallskip\noindent
{\bf Case 2.} $|ad-bc| \ge 2$. In this case
$$\left(\matrix{k\cr \ell\cr}\right) = \left(\matrix{\alp\cr \beta\cr}\right)
+ p\left(\matrix{a\cr b\cr}\right) +
q\left(\matrix{c\cr d\cr}\right)$$
and necessarily $|L\cap \ZZ^2_+|
=\infty$. Now
$$0=\sum_{r=1}^n c_r z_r^k \oz_r^\ell = \sum_{r=1}^n c_r z_r^{\alp+pa +qc}
\oz_r^{\beta +pb +qd} =
\sum_{r=1}^n \left(c_r z_r^\alp \oz_r^\beta\right)\left(z_r^a
\oz_r^b\right)^p \left(z_r^c\oz_r^d\right)^q$$
for all $p,q\in \ZZ$.

We must have $a+b\ne 0$ or $c+d\ne 0$. For otherwise
$ad-bc=0$. Assume, without loss of generality, that $a+b\ne 0$. Fixing
any $q\in\ZZ$, we have
$$0 = \sum_{r=1}^n \left(c_r z_r^\alp \oz_r^\beta \left(z_r^c\oz_r^d\right)^q
\right)\left(z_r^a \oz_r^b\right)^p
= \sum_{r=1}^n d_r w_r^p$$
for all $p\in\ZZ$, where
$$d_r = c_r z_r^\alp \oz_r^\beta (z_r^c\oz_r^d)^q \ne 0$$
and
$$ w_r = z_r^a \oz_r^b\,.$$
As $|w_r| = |z_r|^{a+b}$ and $a+b\ne 0$ we have
$$\{r\,: \,|z_r|=\lam\} =\{r\,:\, |w_r|=\lam^{a+b}\}\,.$$
Applying Proposition 3.10 gives us our result. \eop

\medskip
Our proof of Theorem 3.5 is a consequence of Theorem 3.8 and
Propositions 3.7 and 3.9.

\medskip\noindent
{\bf Proof of Theorem 3.5.} Assume to the contrary that there exists a
$J\subseteq \ZZ^2_+$ such that for each $p\in \NN$ and
$q\in \{0,1\nek p-1\}$ we have
$$\left| \{k-\ell\,:\, (k,\ell)\in J,\, k-\ell =q\,(\mod\, p)\}\right|
= \infty\,,$$
and yet there exist
distinct nonzero points $z_1\nek z_n$ in $\CC$ and
nonzero values $c_1\nek c_n$ in $\CC$ such that
$$\sum_{r=1}^n c_r z_r^k\oz_r^\ell = 0$$
for all $(k,\ell)\in J$. Set
$$\calB = \{(k,\ell)\,:\, (k,\ell)\in\ZZ^2,\,
\sum_{r=1}^n c_r z_r^k\oz_r^\ell = 0\}\,.$$
Then $J\subseteq \calB$ and from Theorem 3.8 we know that $\calB$ is a
union of a finite number of translates of subgroups of $\ZZ^2$.

Let $\calB^*$ be the union of those translates of subgroups of
$\ZZ^2$ in $\calB$ which have an infinite number of elements in
$\ZZ^2_+$. That is,
$$\calB^*=\union_{s=1}^m L_s$$
where $L_s$ is a translate of a subgroup of $\ZZ^2$, $L_s\subseteq
\calB$, and $|L_s\cap \ZZ^2_+|=\infty$. Note that $\calB^*$ contains
all of the elements of $J$, except perhaps for a finite number of
indices. Thus
$$\left| \{k-\ell\,:\, (k,\ell)\in \calB^*,\, k-\ell =q\,(\mod\, p)\}\right|
= \infty\,,$$
for all $p$ and $q$ as above.

Applying Proposition 3.9 to the subgroups $L_s$ of $\calB^*$ we obtain
$$\sum_{\{r:|z_r|=\lam\}} c_r z_r^k\oz_r^\ell = 0$$
for all $(k,\ell)\in\calB^*$. This implies that for each $\lam\in
\RR_+$ we have
$$\sum_{\{r:|z_r|=\lam\}} c_r z_r^{k-\ell} = 0$$
for all $(k,\ell)\in\calB^*$. Applying Proposition 3.7 to
$$J^*=\{k-\ell\,:\, (k,\ell)\in\calB^*\}$$
a contradiction ensues. \eop

\noindent
{\bf Remark 3.} As in Remark 1, the condition on $J$ given in
Theorem 3.1 is the necessary and sufficient so that for every
finite point set $E$ in $\CC$
$$\dim \Pi_J|_E = |E|$$
where
$$\Pi_J=\span \{\, z^k \oz^\ell\,:\, (k,\ell)\in J\,\}\,.$$

\medskip\noindent
{\bf Remark 4.} If we restrict ourselves to $S$, the unit sphere in $H$, i.e.,
$S=\{z\,:\,z\in H, \|z\|=1\}$, then it readily follows from the
above analysis that $f$ of the form (3.1) is strictly Hermitian positive
definite on $S$ if and only if for each $p\in \NN$ and $q\in \{0,1\nek
p-1\}$
$$\left| \{k-\ell:\, (k,\ell)\in J\,,\, k-\ell =q\, (\mod \,p)\}\right| =
\infty\,.$$
In other words, the condition is exactly the same as for all of $H$ except
that we do not have or need $(0,0)\in J$.

\Acknowledgements{We would like to thank Jean-Paul Bezivin, Michel
Laurent, Jozef Siciak, Vilmos Totik and Bronislaw Wajnryb for
their kind assistance and patience.}

\References

\ref Aronszajn, N., Theory of reproducing kernels, {\sl
Trans.~Amer.~Math.~Soc.} {\bf 68} (1950), 337--404.

\ref Berg, C., J.~P.~R.~Christensen, and P.~Ressel, {\sl Harmonic
Analysis on Semigroups}, Springer-Verlag, New York, 1984.

\ref Bernstein, S., Sur l'ordre de la meilleure approximation des
fonctions continues par les polyn\^omes de degr\'e donn\'e, {\sl
Mem.~Cl.~Sci.~Acad.~Roy.~Belg.} {\bf 4} (1912--1922), Fasc.~1
(1912), 1--103.

\ref Bochner, S., Montone Funktionen, Stieltjes Integrale und
harmonische Analyse, {\sl Math. Annalen} {\bf 108} (1933),
378--410.

\ref Browder, F.~E., Real analytic functions on product spaces and
separate analyticity, {\sl Can.~J.~Math.} {\bf 13} (1961), 650--656.

\ref Cameron, R.~H., and D.~A.~Storvick, Analytic continuation for
functions of several complex variables, {\sl Trans.~Amer.~Math.~Soc.}
{\bf 125} (1966), 7--12.

\ref Chen, D., Menegatto, V.~A., and X.~Sun, A necessary and
sufficient condition for strictly positive definite functions on
spheres, {\sl Proc.~Amer.~Math.~Soc.} {\bf 131} (2003), 2733--2740.

\ref Christensen, J.~P.~R., and P.~Ressel, Functions operating on
positive definite matrices and a theorem of Schoenberg, {\sl
Trans.~Amer.~Math.~Soc.} {\bf 243} (1978), 89--95.

\ref Chung, K.~L., {\sl A Course in Probability Theory}, Harcourt,
Brace \& World, New York, 1968.

\ref Donoghue, W. F., {\sl Monotone Matrix Functions and Analytic
Continuation}, Springer-Verlag, New York, 1974.

\ref FitzGerald, C.~H., C.~A.~Micchelli, and A.~Pinkus,
Functions that preserve families of positive semidefinite
matrices, {\sl Linear Alg.~and Appl.} {\bf 221} (1995), 83--102.

\ref Gneiting, T., Criteria of P\'olya type for radial positive
definite functions, {\sl Proc.~Amer. Math.~Soc.} {\bf 129} (2001),
2309--2318.

\ref Herz, C.~S., Fonctions op\'erant sur les fonctions
d\'efinies-positive, {\sl Ann.~Inst.~Fourier (Grenoble)} {\bf 13}
(1963), 161--180.

\ref Horn, R.~A., The theory of infinitely divisible matrices and
kernels, {\sl Trans.~Amer.~Math. Soc.} {\bf 136} (1969), 269--286.

\ref Horn, R.~A.nd C.~R.~Johnson, {\sl Topics in Matrix
Analysis}, Cambridge Univ.~Press, Cambridge, 1991.

\ref Kawata, T., {\sl Fourier Analysis in Probability Theory},
Academic Press, New York, 1972.

\ref Laurent, M., Equations diophantiennes exponentielles, {\sl
Invent.~Math.} {\bf 78} (1984), 299--327.

\ref Laurent, M., Equations exponentielles-polyn\^omes et suites
r\'ecurrentes lin\'eaires. II, {\sl J.~Number Theory} {\bf 31} (1989),
24--53.

\ref Lu, F., and H.~Sun, Positive definite dot product kernels in
learning theory, to appear in {\sl Adv.~Comp.~Math.}

\ref Mercer, J., Functions of positive and negative type, and their
connection with the theory of integral equations, {\sl
Philos.~Trans.~Royal Soc.~London}, Ser.~A, {\bf 209} (1909), 415--446.

\ref Meschkowski, H., {\sl Hilbertsche R\"aume mit Kernfunktion},
Springer-Verlag, Berlin, 1962.

\ref Nguyen, T.~V., Separate analyticity and related subjects,
{\sl Vietnam J.~Math.} {\bf 25} (1997), 81--90.

\ref Pinkus, A., Strictly positive definite functions on a real inner
product space, to appear in {\sl Adv.~Comp.~Math.}

\ref P\'olya, G., and G.~Szeg\"o, {\sl Problems and Theorems in
Analysis, Volume II}, Springer-Verlag, Berlin, 1976.

\ref Rudin, W., Positive definite sequences and absolutely monotonic
functions, {\sl Duke Math.~J.} {\bf 26} (1959), 617--622.

\ref Schoenberg, I.~J., Metric spaces and completely monotone
functions, {\sl Ann.~of Math.} {\bf 39} (1938), 811--841.

\ref Schoenberg, I.~J., Positive definite functions on spheres,
{\sl Duke Math.~J.} {\bf 9} (1942), 96--108.

\ref Shorey, T.~N., and R.~Tijdeman, {\sl Exponential Diophantine
Equations}, Cambridge Tracts in Mathematics, Vol.~87, Cambridge, 1986.

\ref Siciak, J., Separately analytic functions and envelopes of
holomorphy of some lower dimensional subsets of $\CC^n$,
{\sl Ann.~Polon.~Math.} {\bf 22} (1969), 145--171.

\ref Sun, X., and V.~A.~Menegatto, Strictly positive definite functions on the complex Hilbert sphere,
{\sl Adv.~Comp.~Math.} {\bf 11} (1999), 105--119.

\Address
Allan Pinkus
Department of Mathematics
Technion, 32000
Haifa, Israel
pinkus@tx.technion.ac.il

\end